\documentclass[preprint,11pt]{elsarticle}

\usepackage{graphicx}
\usepackage{amsmath,amssymb}
\newcommand{\transv}{\mathrel{\text{\tpitchfork}}}
\makeatletter
\newcommand{\tpitchfork}{%
  \vbox{
    \baselineskip\z@skip
    \lineskip-.52ex
    \lineskiplimit\maxdimen
    \m@th
    \ialign{##\crcr\hidewidth\smash{$-$}\hidewidth\crcr$\pitchfork$\crcr}
  }%
}
\makeatother

 \usepackage{amsthm}
\usepackage[all]{xy}
\usepackage{lineno}
\usepackage{enumerate}
\newtheorem{main}{Theorem}

\newtheorem{thm}{Theorem}[section]

\newtheorem{cor}[thm]{Corollary}
\newtheorem{lem}[thm]{Lemma}
\newtheorem{prop}[thm]{Proposition}
\newtheorem{defn}[thm]{Definition}
\newtheorem{rem}{Remark}[section]

\newtheorem{asmp}{Assumption}
\newcommand{\eps}{\varepsilon}
\newcommand{\chro}{\stackrel{\longrightarrow}{{\rm exp}}\int}

\newcommand{\ad}{\mbox{ad}}

\journal{Journal Name}

\begin{document}

\begin{frontmatter}


\title{Control in the spaces of ensembles of points}


 \author[1adr]{Andrei Agrachev}
\author[2adr]{Andrey Sarychev}
\address[1adr]{Scuola Internazionale degli Studi Avanzati (SISSA), via Bonomea 265,  Trieste, 34136, Italy (agrachev@sissa.it)}
\address[2adr]{Department of Mathematics and Informatics U.Dini, University of Florence, via delle Pandette, 9, Firenze, 50127, Italy (asarychev@unifi.it)}



\begin{abstract}
We study  the controlled dynamics of  the   {\it ensembles of points} of a Riemannian manifold $M$.  Parameterized ensemble of points of $M$ is  the image of a continuous   map   $\gamma:\Theta \to M$, where $\Theta$ is a compact set of parameters. The  dynamics of  ensembles is defined by
the action $\gamma(\theta) \mapsto P_t(\gamma(\theta))$ of the semigroup of diffeomorphisms
$P_t:M \to M, \ t \in \mathbb{R}$, generated by the controlled equation  $\dot{x}=f(x,u(t))$
on $M$.
Therefore any control system  on $M$ defines  a control system on (generally  infinite-dimensional)  space $\mathcal{E}_\Theta(M)$ of the ensembles of points.

We wish to establish criteria of  controllability for  such   control systems.  As in our previous work (\cite{ABS})  we seek to adapt
the   Lie-algebraic approach of geometric control theory  to the infinite-dimensional setting. We  study the case of finite ensembles and prove  genericity of exact controllability property for them.
We  also  find sufficient   approximate controllability criterion for continual ensembles and  prove a result on motion planning in the space  of flows on $M$.
We discuss the relation of the obtained controllability criteria  to various  versions of Rashevsky-Chow theorem for finite- and infinite-dimensional manifolds.

\end{abstract}


\end{frontmatter}



\section{Introduction and problem setting}

 Let $M$ be 
 $C^\infty$-smooth $n$-dimensional ($n \geq 2$) connected Riemannian manifold, with $d(\cdot , \cdot)$, being the Riemannian distance. Let  $\mathcal{E}_\Theta(M)$ be the space of  continuous maps  $\gamma:\Theta \to M$, where $\Theta$ is a compact Lebesgue measure set.
    We call the elements of $\mathcal{E}_\Theta(M)$ {\it ensembles of points} or, for  brevity,  {\it ensembles}.
The space $\mathcal{E} _\Theta(M)$
is infinite-dimensional, whenever       $\Theta$ is an infinite set (see Section \ref{baman}).



     In the control-theoretic setting one looks at  the action on $\mathcal{E}_\Theta(M)$
     of the group  of diffeomorphisms of $M$,
     which are generated by the vector fields from the family   $\{f^u| \ u \in U\} \subset \mbox{\rm Vect }M $. Alternatively we can consider the action of  the flows,  defined by the  controlled  equations
        \begin{equation}\label{ncf}
      \dot{x}=f(x,u(t)), \ u(t) \in U,
         \end{equation}
where $u(t)$ are admissible, for example,  piecewise-constant, or piecewise-continuous, or boundary measurable controls.

      The flow $P^{u(\cdot)}_t \ (P_0=Id)$, generated by  control system (\ref{ncf}) and a  given admissible control $u(t)=(u_1(t), \ldots , u_r(t))$, acts  on $\gamma (\theta) \in  \mathcal{E}_\Theta(M)$ according to the formula
  \[\hat{P}^{u(\cdot)}_t : \gamma(\theta) \mapsto P^{u(\cdot)}_t(\gamma (\theta)), \theta \in \Theta.\]
Thus control system (\ref{ncf}) gives rise to  a  control system in the space of ensembles $\mathcal{E} _\Theta(M)$.
   We set the controllability problem for the action of   control system (\ref{ncf}) on  $\mathcal{E}_\Theta(M)$.

      \begin{defn}\label{codef}
Ensemble $\alpha(\cdot) \in \mathcal{E}_\Theta(M)$ can be steered  in time-$T$
 to ensemble $\omega(\cdot) \in \mathcal{E}_\Theta(M)$  by  control system (\ref{ncf}), if there exists a control
$\bar{u}\in L_\infty([0,T],U)$   such that for the flow $P^{\bar{u}(\cdot)}_t$, generated by the equation $\dot{x}=f(x(t),\bar{u}(t))$,
there holds
 \[ P^{\bar{u}(\cdot)}_T\left(\alpha(\theta)\right)=\omega(\theta). \ \Box \]
\end{defn}

\begin{defn}\label{atset} The time-$T$ attainable set from $\alpha(\cdot) \in \mathcal{E}_\Theta(M)$ for  control system (\ref{ncf}) in the space of ensembles $\mathcal{E}_\Theta(M)$ is   \[\mathcal{A}_T(\alpha(\cdot))=\{P^{u(\cdot)}_T\left(\alpha(\theta)\right)| \ u(\cdot) \in  L_\infty([0,T],U)\} \subset \mathcal{E}_\Theta(M).\ \Box\]
\end{defn}

\begin{defn}\label{glocon}
Control system (\ref{ncf}) is globally exactly controllable in time $T$  in the space $\mathcal{E}_\Theta(M)$ from  $\alpha(\theta) \in \mathcal{E}_\Theta(M)$, if $\mathcal{A}_T(\alpha (\theta))=\mathcal{E}_\Theta(M)$.
Control system (\ref{ncf}) is time-$T$ globally exactly controllable if it is  globally exactly controllable in time-$T$ from each $\alpha(\theta) \in \mathcal{E}_\Theta(M) \ \Box$.
\end{defn}

\begin{rem}
If $\Theta=\{\theta\}$ is a singleton, then the time-$T$ attainable sets $\mathcal{A}_T(\alpha_\theta)$ coincide
with the standard attainable sets of  system (\ref{ncf}) from the point $\alpha_\theta \in M$.   The notions of global and global approximate controllability coincide with the standard notions for control system (\ref{ncf}) on $M$.
\end{rem}

If $\Theta$  is an infinite set,   it is hard to achieve {\it exact} ensemble controllability for system \eqref{ncf}.
Instead we will  study  $C^0$- or $L_p$-approximate controllability property.

      \begin{defn}\label{acodef}
Ensemble $\alpha(\cdot) \in \mathcal{E} _\Theta(M)$ is
$C^0$-approximately steerable in time-$T$ to ensemble  $\omega(\cdot) \in \mathcal{E}_\Theta(M)$  by control system (\ref{ncf}), if for each $\eps >0$ there  exists $\bar u(\cdot)$ such that
\begin{equation}\label{c0_ap}
\sup_{\theta \in \Theta} d\left(\omega(\theta) , P^{\bar{u}(\cdot)}_T\left(\alpha(\theta)\right)\right) \leq \eps .
\end{equation}

Ensemble $\alpha(\cdot) \in \mathcal{E} _\Theta(M)$ is
$L_p$-approximately steerable in time-$T$ to ensemble $\omega(\cdot) \in \mathcal{E} _\Theta(M)$  by
 control system (\ref{ncf}), if for each $\eps >0$ there  exists $\bar u(\cdot)$ such that
\[\int_\Theta \left( d\left(\omega(\theta),P^{\bar{u}(\cdot)}_T\left(\alpha(\theta)\right)\right)\right)^p d\theta \leq \eps^p .\ \qed.\]
\end{defn}

\begin{defn}
 Control system (\ref{ncf}) is time-$T$ globally approximately controllable from $\alpha(\cdot) \in \mathcal{E}_\Theta(M)$ if $\mathcal{A}_T(\alpha)$ is dense in
 $\mathcal{E}_\Theta(M)$ in the respective metric. The  system  is time-$T$ globally approximately controllable if it is time-$T$ globally approximately controllable from each $\alpha(\cdot) \in \mathcal{E}_\Theta(M). \ \Box$
\end{defn}

It is known that the attainable sets and the controllability properties of  control system (\ref{ncf}) on $M$
can be characterized via  properties of the  Lie brackets of the vector fields $f^u(x), \ u \in U$. In particular case for  a symmetric control-linear system
      \begin{equation}\label{cls}
     \dot{x}=\sum_{j=1}^sf_j(x)u_j(t)
     \end{equation}
global controllability property {\it for singletons} is guaranteed by  the  {\it bracket generating condition}: for  each point $x \in M$ the evaluations at $x$  of the iterated Lie brackets $[f_{j_1},[ \ldots [f_{j_{N-1}},f_{j_{N}}]\ldots ]$ span the tangent space $T_xM$.

We are going to establish  controllability criteria for control system (\ref{cls}) acting in the space of  ensembles $\mathcal{E}_\Theta(M)$.  The criteria for finite and continual ensembles are provided in  Sections \ref{fipoen} and \ref{rccon}.
As far as controlled dynamics in the space of ensembles is defined by  action of  the flows, generated by
controlled system \eqref{cls}, it is important to analyze whether and how the controllability criterion could be "lifted" to the group of diffeomorphisms
or the semigroup  of flows. This is done  in Section~\ref{ext_flows}, where
 Theorem~\ref{liex_flows} provides a result on a Lie extension of the action of  system \eqref{cls} in the group of diffeomorphisms.
 In Section~\ref{disc_chow} we discuss the relation of the established controllability criterion for continual ensembles of points
to various  versions of  Rashevsky-Chow theorem in finite and infinite dimensions. It turns out that the latter typically are not applicable to ensemble controllability.

The proofs of the main results are provided in Sections   \ref{proof1}-\ref{proof2}.

By now there are numerous   publications on  simultaneous control of ensembles  of control systems
\begin{equation}\label{systheta}
  \dot{x}=f(x,u,\theta), \ x \in M, \ u \in U, \theta \in \Theta
\end{equation}
by a unique control.
This direction of study  has been initiated by S.~Li and N.~Khaneja (\cite{LiKh05,LiKh09}) for the case of quantum ensembles.
Few other  publication which took on the subject are \cite{BCR,BST,FH}, where readers can find more bibliographic references.
In our previous publication \cite{ABS} we considered the ensembles of systems \eqref{systheta}, and formulated Lie algebraic controllability criteria for ensembles of systems.

In the present publication  we consider ensembles of points controlled by virtue of a {\em single system} and {\em single open loop control}.
This choice  distinguishes the  problem setting not only from  the previous one, but also  from the control problems, in  which both the state space and the set of control parameters are infinite-dimensional. Examples of the latter kind appear  in \cite{AgCa} and are common  in the literature on mass transportation.

A common feature which this  publication shares with \cite{ABS} is the Lie algebraic approach to study of controllability; we seek to demonstrate  parallelism in formulations  along the text.

\section{Banach manifold of ensembles}
\label{baman}

As we said ensembles of points in $M$ are the  images of  continuous   maps   $\gamma:\Theta \to M$;
      the set  of parameters $\Theta$  is assumed to be    compact.
At some moments   we  assume additionally  the maps $\gamma$  to be    injective.
The set of  ensembles is denoted by $\mathcal{E}_\Theta(M)$.

Whenever the set of parameters $\Theta$ is  finite,
then the ensemble is called finite and  the set of ensembles $\mathcal{E}_\Theta(M)$ is a finite-dimensional manifold.

Define for any ensemble $\gamma(\theta)  \in \mathcal{E}_\Theta(M)$ a tangent space $T_\gamma\mathcal{E}_\Theta(M)$, consisting of the continuous maps $T\gamma: \Theta \rightarrow TM$, for which the diagram
\[\xymatrix{\Theta \ar[rr]^-{T\gamma} \ar[dr]_-\gamma &&TM \ar[dl]_-\pi  \\
& M &
         }       \]
is commutative. Representing an element of the tangent bundle  $TM$ as a pair $(x,\xi), \ x \in M, \xi \in T_x M$, we note that
\[T\gamma(\theta)=(\gamma(\theta), \xi(\theta)), \ \xi(\theta) \in T_{\gamma(\theta)}M, \ \theta \in \Theta .\]
If $M=\mathbb{R}^n$, then  $T_\gamma\mathcal{E}_\Theta(M)$ can be identified with the set of continuous maps $C^0(\Theta , \mathbb{R}^n \times \mathbb{R}^n)$.

One can define a vector field on $\mathcal{E}_\Theta(M)$ as a section of the  tangent  bundle $T\mathcal{E}_\Theta (M)$.

The flow $e^{tf}$, generated by a time-invariant vector field   $f \in \mbox{Vect}(M)$, and acting onto an ensemble $\gamma(\theta)$,   defines a lift of $f$ to
the vector field $$F \in \mbox{Vect}\left(\mathcal{E}_\Theta(M)\right): \ F(\gamma(\cdot))=\left.\frac{d}{dt}\right|_{t=0}e^{tf}\left(\gamma(\cdot)\right)=f(\gamma(\cdot)).$$
      The same holds for   time-variant vector fields  $f_t$.

   The Lie brackets of the lifted vector fields are the lifts of the Lie brackets of the vector fields: $\left.[F_1,F_2]\right|_{\gamma(\cdot)}=[f_1, f_2](\gamma(\cdot))$.

One can provide $T_{\gamma (\cdot)}\mathcal{E}_\Theta(M)$ with different metrics.   Of interest  for us are those  obtained by the restrictions
of the  metrics $C^0(\Theta , TM)$, and $L_p(\Theta , TM)$ onto $T\mathcal{E}_\Theta(M)$.

\section{Genericity of the controllability property  for finite ensembles of points}
\label{fipoen}

 Let $\Theta=\{1, \ldots , N \}$. Finite ensemble $\gamma : \Theta \mapsto M$  is  an $N$-ple of  points $\gamma =(\gamma_1, \ldots , \gamma_N) \in M^N$.
In this Section  we assume $\gamma$ to be injective, so that the points $\gamma_j$ are pairwise distinct.
Let $\Delta^N \subset M^N$ be the set of $N$-ples $(x_1, \ldots , x_N) \in M^N$ with (at least) two coinciding components:   $x_i=x_j$, for some $i \neq j$.
Then the space of ensembles $\mathcal{E}_N(M)$ is  identified  with the complement of $\Delta^N: \ \mathcal{E}_N(M)=M^N \setminus \Delta=M^{(N)}$.

For each $\gamma \in M^{(N)}$ the tangent space $T_\gamma M^{(N)}$ is isomorphic to
\[\bigotimes_{j=1}^NT_{\gamma_j}M =T_{\gamma_1}M \times \cdots \times T_{\gamma_N}M.\]

For a vector field $X \in \mbox{\rm Vect} M$  consider its $N$-fold, defined on $M^{(N)}$  as   $X^N(x_1, \ldots , x_N)=(X(x_1), \ldots ,X(x_n))$. For $X,Y \in \mbox{\rm Vect } M$, and $N \geq 1$ we define the Lie bracket of the $N$-folds $X^N,Y^N$  on $M^{(N)}$  "componentwise":
$[X^N,Y^N]=[X,Y]^N$, where $[X,Y]$ is the  Lie bracket of $X,Y$ on $M$.
The same holds for the iterated Lie brackets.

Given the vector fields $f_1, \ldots f_s$ on $M$ their $N$-folds $f_1^N, \ldots , f_s^N$  form a
 bracket generating system on $M^{(N)}$, if the evaluations
 of their iterated Lie brackets  at each $\gamma \in M^{(N)}$, span the tangent space $T_\gamma M^{(N)}=\bigotimes_{j=1}^NT_{\gamma_j}M$.
Evidently  for $N>1$ the property strictly is stronger, than the bracket generating property for $f_1, \ldots , f_s$ on $M$.
We provide some comments below in Section \ref{disc_chow}.

The following result is a  corollary of classical Rashevsky-Chow theorem (see Proposition \ref{rct}).


\begin{prop}[global controllability criterion for system \eqref{cls} in the space of  finite point ensembles]
\label{glec}
If the $N$-folds $f_1^N, \ldots , f_s^N$ are bracket generating at each point of $M^{(N)}$,
then $\forall T>0$ the system (\ref{cls}) is time-$T$ globally exactly controllable in the space of finite
ensembles $(\gamma_1, \ldots , \gamma_N) \in M^{(N)}. \ \Box$
\end{prop}


 Proposition \ref{glec} relates global controllability of  system (\ref{cls}) for $N$-point ensembles
 to the  bracket generating property on $M^{(N)}$ for the $N$-folds  of the vector fields $f_1, \ldots , f_s$.
  The following result states that the bracket generating  property for $N$-folds is generic.

\begin{main}\label{generic}
  For any $N \geq 1$ and  sufficiently  large $\ell$,  there is  a set of $s$-ples of vector fields $(f_1, \ldots , f_s)$, which is residual  in 
  $\mbox{\rm Vect} M ^{\otimes s}$
  in Whitney $C^\ell$-topology,   such that for any
  $(f_1, \ldots , f_s)$ from this set  the
  $N$-folds $(f_1^{N}, \ldots , f_s^{N})$
  are bracket generating at each point of $M^{(N)}=M^N \setminus \Delta^N$. $\qed$
   \end{main}

Note that the notion of genericity in the theorem allows for (small) perturbations of the  $f_i$, but not of $f_i^N=(f_i, \ldots , f_i)$  directly.
Therefore the theorem extends the
classical result by  C.Lobry's (\cite{Lob}) on genericity of the property of  controllability for singletons (see also Theorem 3.1 of our previous work  \cite{ABS}
 on the genericity of controllability property for  ensembles  of control systems).

Proof of theorem \ref{generic}  (for $s=2$) is  provided in Section \ref{proof1}.


\section{Criterion of approximate steering  for continual ensembles of points}
\label{rccon}

 To formulate  criterion  for approximate steering of  continual ensembles of points we impose  the following  assumption for control system \eqref{cls}.
\begin{asmp}[boundedness in $x$]
\label{unal}
      The $C^\infty$-smooth vector fields   $f_j(x) \in \mbox{Vect }M, \ j=1, \ldots , s$,  which define system \eqref{cls},  are bounded on  $M$
      together with their covariant
      derivatives of each order.  \hfill $\qed$
\end{asmp}

The boundedness of $f_j$  and of their covariant derivatives on $M$ implies  completeness of the vector fields $f_j$ and of their Lie brackets of any order.
Completeness of a vector field means that
the  trajectory of  the vector field  with arbitrary  initial  data  can be extended to  each  compact subinterval of the time axis.

This assumption is rather natural. It holds for  compact manifolds $M$. For a non-compact $M$ it obviously holds for vector fields with compact supports.
Other  examples are vector fields on $\mathbb{R}^n$, whose components are trigonometric polynomials in $x$, or polynomial (in $x$) vector fields, multiplied by  functions rapidly decaying at infinity (e.g. by $e^{-x^2}$).

Consider a couple of  initial and target  ensembles of points $\alpha(\theta), \omega(\theta) \in \mathcal{E}_\Theta(M)$, which we assume to be
diffeotopic,\footnote{We can assume instead an existence,   for each $\eps >0$,   of an ensemble  $\omega_\varepsilon (\cdot)$, which is $\eps$-close to $\omega(\cdot)$ in $C^0(\Theta)$-metric and
    diffeotopic to $\alpha(\cdot)$.} i.e. satisfying the relation $R_T(\alpha(\cdot))=\omega(\cdot)$,
 where $t \to R_t, \ t \in [0,T], \ R_0=\mbox{Id}$,  is a flow on $M$,  defined by a time-variant vector field $Y_t(x)$, with $Y_t(x), D_xY_t(x)$ continuous.

Note that  the (reference) flow $R_t$ is a priori unrelated to  control system (\ref{cls}).
Denote by $\gamma_t(\theta)$ the image of $\alpha(\theta)$ under the diffeotopy \[ \gamma_t(\theta)=R_t(\alpha(\theta)), \ \gamma_0(\theta)=\alpha(\theta), \ \gamma_T(\theta)=\omega(\theta).\]

 We introduce standard notation for the seminorms in the space of vector fields on $M$: for a compact $K \subset M$
 \[ \|X\|_{r,K} =\sup_{x \in K}\left(\sum_{0 \leq |\beta| \leq r }\left|D^\beta X(x)\right|\right)  \]
and
 \[ \|X\|_{r} =\sup_{x \in M}\left(\sum_{0 \leq |\beta | \leq r }\left|D^\beta X(x)\right|\right) .  \]

 Let $\mbox{Lie}{\{f\}}$
 be  the Lie algebra,  generated   by the vector fields $f_1, \ldots , f_s$.  Put  for $\lambda>0$ and a compact $K \subset M$:
  \[\mbox{Lie}^\lambda_{1,K}\{f\}=\left\{X(x) \in \mbox{Lie}\{f\}\left| \
  \|X\|_{1,K} < \lambda \right.\right\},\]
  and
  \[ \mbox{Lie}^\lambda_{1}{\{f\}}=\left\{X(x) \in \mbox{Lie}{\{f\}}\left| \
  \|X\|_{1} < \lambda \right.\right\}.\]


The following {\it bracket approximating  condition along  a  diffeotopy} is  the key part
of the criterion for steering continual ensembles of points.   In Section~\ref{disc_chow} we discuss the
reason for the choice of this particular form of condition.

\begin{defn}[Lie bracket $C^0$-approximating condition along a diffeotopy]\label{def_brapr}
 Let  the diffeotopy $\gamma_t=R_t(\alpha(\cdot)), \ t \in [0,T]$,
generated by the vector field  $Y_t(x)$,  join $\alpha(\cdot)$ and $\omega(\cdot)$.   System  (\ref{cls}) satisfies  Lie bracket $C^0$-approximating  condition along $\gamma_t$,    if there exist $\lambda>0$ and   a compact neighborhood $\mathcal{O}_\Gamma$ of the set
$\Gamma=\{\gamma_t(\theta)| \ \theta \in \Theta , \ t \in [0,T] \}$  such that
\begin{equation}\label{C0_appr}
\forall t \in [0,T]: \   \inf \left\{\left. \sup_{\theta \in \Theta}\left|Y_t(\gamma_t(\theta))-X(\gamma_t(\theta)) \right| \  \right|
\ X \in  \mbox{Lie}^\lambda_{1,\mathcal{O}_\Gamma}{\{f\}}\right\}=0.
\end{equation}
\end{defn}

\begin{main}[approximate steering  criterion for   ensembles of points]\label{trc} Let $\alpha(\theta), \omega(\theta)$
be two  ensembles of points, joined by a diffeotopy $\gamma_t(\theta), \ t \in [0,T]$. If  control system (\ref{cls}) satisfies the
Lie bracket $C^0$-approximating condition along the   diffeotopy,   then  $\alpha (\cdot)$
      can be  steered $C^0$-ap\-pro\-ximately to $\omega(\cdot)$    by   system  (\ref{cls})   in time $T. \  \square$
\end{main}

\subsection{Approximate controllability for continual ensembles: basic example}
\label{basic_example}

We provide an  example of  application of  Theorem \ref{trc}.
 Consider  system in $\mathbb{R}^2$ with two controls:
\begin{equation}\label{e_R2_system}
  \dot{x_1}=u, \   \dot{x_2}=\varphi(x_1)v.
\end{equation}
        It is a particular case of the control-linear  system  \eqref{cls}:
         \begin{equation}\label{clin}
           \dot{x}=f_1(x)u+f_2(x)v, \ f_1=\partial/\partial x_1, \ f_2=\varphi(x_1)\partial/\partial x_2 .
         \end{equation}
         We assume  $\varphi(x_1)$ to be $C^\infty$-smooth. In our example $\varphi(x_1)=e^{-x_1^2}$.

Choose the initial ensemble
\begin{equation}\label{example_ex2}
  \alpha(\theta)=(\theta , 0), \
 \theta \in \Theta=[0,1].
 \end{equation}

If one takes  for example $u=0$  in \eqref{e_R2_system}, then  $x_1$ remains fixed,  and  by \eqref{e_R2_system},\eqref{example_ex2}
$$x_2(T;\theta)=m_{v(\cdot)}\varphi(\theta),$$
where  $m_{v(\cdot)}=\int_0^T v(t)dt \ \in  \mathbb{R}$.
Therefore for vanishing $u(\cdot)$ the set of "attainable profiles" for $x_2(T;\theta)$ is very limited.

To illustrate Theorem \ref{trc} we fix target ensemble $\omega(\theta)=(\theta , \theta)$ and
 choose a  diffeotopy
\begin{equation}\label{diff_gam}
  \gamma_t(\theta)=(\theta,  t \theta), \ t \in [0,1],
\end{equation}
which joins $\alpha(\theta)$ and $\omega(\theta)$.
The diffeotopy  is generated by the (time-invariant) vector field  $Y(x)=Y(x_1,x_2)=x_1 \partial/\partial x_2$.
Evaluation of the  vector field $Y$  along the diffeotopy \eqref{diff_gam},
equals
$$\forall t \in [0,1]: \ Y(\gamma_t(\theta)=Y(\theta,  t \theta)=\theta \partial/\partial x_2.$$

                  The Lie algebra, generated by $f_1,f_2$,  is spanned in the treated  case by the vector fields
                  \begin{equation}\label{lie_alg}
                   f_1, \ad^k f_1 f_2 =\varphi^{(k)}(x_1)\partial/\partial x_2, \ k=0,1, 2 , \ldots
                  \end{equation}
and is  infinite-dimensional for our choice of $\varphi(\cdot)$.


The evaluations $ f_1(\gamma_t(\theta))$ and $\ad^k f_1 f_2(\gamma_t(\theta))$
equal
$$f_1(\gamma_t(\theta))=\partial/\partial x_1, \  \left(\ad^k f_1 f_2\right)(\gamma_t(\theta))=\varphi^{(k)}(\theta)\partial/\partial x_2, \ k=0,1,2 \ldots  .$$

The  successive derivatives of $\varphi(x)=e^{-x^2}$ are
 \begin{equation}\label{herm}
 \varphi^{(m)}(x)=(-1)^mH_m(x)e^{-x^2}, \ m=0,1, \ldots  ,
 \end{equation}
   where $H_m(x)$ are Hermite polynomials. Recall that $H_m(x)$  form an orthogonal complete system  for   $L_2(-\infty, +\infty)$ with the weight $e^{-x^2}$.

Let $\mathcal{H}$ be  (infinite-dimensional) linear space generated by functions \eqref{herm}. Generic element of  $\mbox{Lie}\{f_1, f_2\}$  can be represented as
$$a \frac{\partial}{\partial x_1}+h(x_1) \frac{\partial}{\partial x_2}, \ a \in \mathbb{R}, \ h \in \mathcal{H} ,$$
and its evaluation at $\gamma_t(\theta)$ equals $$a \frac{\partial}{\partial x_1}+h(\theta) \frac{\partial}{\partial x_2}, \ a \in \mathbb{R}, \ h \in \mathcal{H} .$$

The $C^0$ bracket approximating condition along $\gamma_t(\theta)$  amounts to
the approximability  in $C^0[0,1]$  of the function $Y_2(\theta)=\theta$ by  the functions from
a bounded equi-Lipschitzian subset of $\mathcal{H}$.

To establish approximability   for  chosen  example   we use the following  technical lemmae.
\begin{lem}\label{equilip}
There exists $\lambda >0$  such that
\[\inf \left\{\left.  \sup_{\theta \in [0,1]} \left|\theta-h(\theta)\right| \  \right|
\ h(\cdot) \in \mathcal{H}, \ \sup_{\theta \in [0,1]}\left( |h(\theta)|+|h'(\theta)|\right)  < \lambda   \right\}=0. \]
\end{lem}

The lemma  follows from  standard facts, which  concern  the  expansions with respect to  Hermite system.

\begin{lem}\label{hermprop}
Let $g(x)$ be a smooth function with compact support in $(-\infty, +\infty)$ and
\begin{equation}\label{hexp}
 g(x) \simeq \sum_{m \geq 0} g_mH_m(x)
\end{equation}
be its expansion with respect to Hermite system.
Then:
\begin{enumerate}[(i)]
  \item expansion \eqref{hexp} converges to $g(x)$ uniformly on any compact interval;
  \item the expansion $\sum_{m \geq 1}g_mH'_m(x)$ converges to $g'(x)$  uniformly on any compact interval.
\end{enumerate}
\end{lem}
For (i) see e.g. \cite[\S 8]{NUv}.  Statement (ii) follows easily from (i),  given the   relation
$  H'_m(x)=2m H_{m-1}(x)$ for the Hermite polynomials.
 Indeed  $$\sum_{m \geq 1}g_mH'_m(x)=\sum_{m \geq 1}2mg_mH_{m-1}(x)=\sum_{m \geq 0}2(m+1)g_{m+1}H_{m}(x),$$
 and it rests to verify that the coefficients of the expansion of $g'(x)$ with respect to  Hermite system are precisely $2(m+1)g_{m+1}$.  This in its turn  follows by
  direct computation by the formulae
$$ \ H'_m(x)=2xH_m(x)-H_{m+1}(x), \  \int_{-\infty}^{+\infty}(H_m(x))^2e^{-x^2}dx=2^m m!\sqrt{\pi}.$$

Now in order to prove  Lemma \ref{equilip} we take a $C^\infty$ smooth function $g(\theta)$ with compact support on $ (-\infty, +\infty)$, whose  restriction to $[0,1]$ coincides
with   the function $y(\theta)=\theta e^{\theta^2}$.
By Lemma \ref{hermprop} (i)  the  expansion $g(\theta) \simeq \sum_mg_mH_m(\theta)$
 converges uniformly on $[0,1]$ to $\theta e^{\theta^2}$,  and hence the series
$ \sum_mg_mH_m(\theta) e^{-\theta^2}$ converges
to $\theta$  uniformly on $[0,1]$.

Differentiating  $\sum_m c_mH_m(\theta)e^{-\theta^2}$ termwise in $\theta$ we get
\[\sum_m c_mH'_m(\theta)e^{-\theta^2} -\sum_m c_mH_m(\theta)2\theta  e^{-\theta^2}. \]
By Lemma \ref{hermprop} (i) and (ii)      the series $\sum_{m \geq 1} c_mH'_m$  and $\sum_{m \geq 0}  c_mH_m(\theta)$ converge uniformly
on $[0,1]$ to bounded functions; the partial sums of these series are equibounded and therefore partial sums of the series
$ \sum_mg_mH_m(\theta) e^{-\theta^2}$ are equiLipschitzian, what concludes the proof of  Lemma \ref{equilip}.

\section{Lie extensions  and approximate controllability for  flows}
\label{ext_flows}

The proof of Theorem~\ref{trc},   provided in   Section  \ref{proof2},
is based on  an infinite-dimensional version of the  method of Lie extensions (\cite{Jurdj, ABS,AgSach}).

According to this method one
 starts with  establishing
the  property of $C^0$-approximate steering by means of  an  extended control fed into    an extended (in comparison with \eqref{cls})  control system
   \begin{equation}\label{ext_ens}
\frac{d x(t)}{dt}=
\sum_{\beta \in B}X^\beta(x)v_\beta(t),
\end{equation}
   where     $X^\beta(x)$ are the iterated Lie brackets
\begin{equation}\label{bralf}
X^\beta(x)= [f_{\beta_1},[f_{\beta_2},[ \ldots ,f_{\beta_N}]\ldots ]](x)
\end{equation}
of the vector fields  $f_1, \ldots , f_s$ (we assume by default,  that  the  vector fields $f_j(x)$   are included into the  family $\{X^\beta(x), \ \beta \in B\}.$) In \eqref{ext_ens}-\eqref{bralf}  the multiindices $\beta =(\beta_1, \ldots , \beta_N)$  belong to a finite subset   $B \subset \bigcup_{N \geq 1} \{1, \ldots , s\}^N$,  and
   $(v_\beta(t))_{\beta \in B}$  is a (high-dimensional) {\it extended control}.

 After the first step  one has to prove  that  the action of the flow, generated by  extended system \eqref{ext_ens}  on $\mathcal{E}_\Theta(M)$ ,  can be approximated by the action of the flow of  system (\ref{cls}), driven by a low-dimensional  control $u(\cdot)=(u_1(\cdot), \ldots , u_s(\cdot))$.
The latter step   is  the core of the method of Lie extensions.

 Therefore  it is reasonable and  desirable to establish the "lifted" approximate controllability result  for  flows on $M$. Such result  not only will allow  to prove  theorem \ref{trc},  but  is certainly interesting on its own.

The respective formulation is given by
\begin{main}\label{liex_flows}
  Let $P^{v(\cdot)}_t$ be   a flow on $M$, generated by extended control system \eqref{ext_ens} and an
   extended  control $v(t)=\left(v_\beta(t)\right)_{\beta \in B}, \ t \in [0,T]$.
  For each $\eps >0$, $r \geq 0$  and   compact $K \subset M$
  there exists  an appropriate  control $u(t)=(u_1(t), \ldots , u_s(t))$ such that the flow $P^{u(\cdot)}_t$,    generated
by  control system \eqref{cls} and the control $u(\cdot)$,    satisfies:
\[ \|P^{v(\cdot)}_t -P^{u(\cdot)}_t\|_{r,K} < \eps , \ \forall t \in [0,T].       \hspace*{2cm}          \]
\end{main}

An obvious application of this theorem to the case of ensembles provides  the following

\begin{cor}
\label{liex_ens}
If  the   ensemble  $\alpha(\theta)$    can be steered approximately to the ensemble $\omega(\theta)$     in time $T$  by an  extended
      system  \eqref{ext_ens}, then the same
can be accomplished by  the original control system \eqref{cls}. $\qed$
\end{cor}

Indeed let $v(\cdot)$ be an extended control for  extended system  \eqref{ext_ens}, such that for the corresponding  flow
$P^{v(\cdot)}_t$ we get $\sup_{\theta \in \Theta}d\left(\omega(\theta),P^{v(\cdot)}_T(\alpha(\theta))\right) < \eps /2$.
By theorem \ref{liex_flows} there exists a control $u(\cdot)$ for system \eqref{cls} such that
$\sup_{\theta \in \Theta}d\left(P^{v(\cdot)}_T(\alpha(\theta)), P^{u(\cdot)}_T(\alpha(\theta))\right) < \eps/2$ and hence
   $$\sup_{\theta \in \Theta}d\left(\omega(\theta),P^{u(\cdot)}_T(\alpha(\theta))\right) < \eps .$$


\section{Theorem \ref{trc} and Rashevsky-Chow theorem(s):  discussion of the formulations}
\label{disc_chow}
The formulations of  the  results, provided in the  two previous sections,  show    similarity to the
formulations of Rashevsky-Chow theorem  on finite-dimensional and infinite-dimensional manifolds.
In this Section we survey these  formulations and establish their relation to Theorem \ref{trc}.

\subsection{Lie rank/bracket generating controllability criteria}
\label{formulations_RC}

Rashevsky-Chow theorem   provides sufficient (and necessary in real analytic case) criterion for global exact  controllability of  system (\ref{cls})
for singletons (= single-point ensembles)  on a connected finite-dimensional manifold $M$   in terms of {\it bracket generating property}.
This property holds for control system (\ref{cls}) at $x \in M$ if  the evaluations of the iterated Lie brackets \eqref{bralf}
of the vector fields  $f_1, \ldots , f_r$ at $x$ span the respective tangent space $T_xM$.

\begin{prop}[Rashevsky-Chow theorem in finite dimension, \cite{AgSach},\cite{Jurdj}]\label{rct}
  Let for control system (\ref{cls})  the  bracket generating property  hold  at each point of $M$.
      Then  $ \forall x_\alpha,x_\omega \in M, \ \forall T>0$  the point $x_\alpha$ can be connected with $x_\omega$
      by an admissible trajectory $x(t), \ t \in [0,T]$ of  system (\ref{cls}), i.e.  system (\ref{cls}) is globally controllable in any time $T$. If the manifold $M$ and the vector fields $f_1, \ldots , f_s$ are real analytic then the bracket generating property is necessary and sufficient for  global controllability of  system (\ref{cls}). $\square$
\end{prop}


The bracket generating property for $f_1, \ldots , f_s$ is by no means sufficient for controllability
of ensembles, even  finite ones.     For example if  this property holds but the Lie algebra $\mbox{Lie}\{f\}$,  correspondent
to  the system (\ref{cls})  is finite-dimensional, then  the $N$-fold of  system  (\ref{cls}) can not possess bracket  generating property
on $M^{(N)}$ (see Section \ref{fipoen}),  if   $N   \dim \ M > \dim \  \mbox{Lie}\{f\}$.  Hence   if $\dim \mbox{Lie}\{f\} < + \infty$,
 then  exact controllability in the space of $N$-point ensembles,  with $N$ sufficiently large, is not achievable.

    What for  continual ensembles, they form, as we said,  infinite-dimensional
    Banach manifold $\mathcal{E}_\Theta(M)$  (see Sections~\ref{baman} and \ref{rccon}) and
    control system   (\ref{cls}) admits a lift to   a control system  on  $\mathcal{E}_\Theta(M)$.

 One can think of application of  infinite-dimensional  Rashevsky-Chow theorem (\cite{DuSam},\cite{Led2}) to the lifted system.

\begin{prop}[infinite-dimensional version of Rashevsky-Chow theorem]
 Consider a control system
     $\dot{y}=\sum_{j=1}^s F_j(x)u_j(t)$,
     defined on Banach manifold $\mathcal{E}$. If  the condition
  \begin{equation}\label{bragen}
    \overline{Lie\{F_1, F_2,\ldots , F_m\}(y)} = T_y \mathcal{E}, \forall y \in \mathcal{E}
  \end{equation}
 holds, then   this system is globally approximately controllable, i.e. for each starting point
 $\tilde y$ the set of points, attainable from $\tilde y$   (by virtue  of the system) is dense in $\mathcal{E}. \ \Box$.    \end{prop}

 Seeking to apply this   result  to the  case of ensembles $\mathcal{E}=\mathcal{E}_\Theta(M)$ one meets two difficulties.

 First,  verification of
 the (approximate)  bracket generating property \eqref{bragen}  has to be done for   {\it  each}  $\gamma(\cdot) \in \mathcal{E}_\Theta(M)$
 and this
results   in a vast set of conditions,  "indexed" by the elements of the functional space $\mathcal{E}_\Theta(M)$.

This difficulty can be overcome by passing to a  pathwise version of Rashevsky-Chow theorem, which in the case of singletons
is  close to its classical formulation.

\begin{prop}\label{prct} Let   $M$  be  a finite-dimensional manifold, $x_\alpha , x_\omega \in M$.
If bracket generating property holds at each  point of  a continuous path $\gamma(\cdot)$, joining   $x_\alpha$ and $x_\omega$,  then $x_\alpha$  and $ x_\omega$   can be joined by an admissible trajectory of $(\ref{cls})$. $\square$
\end{prop}

This result can be deduced directly from Proposition~\ref{rct}.
Indeed if the bracket generating property holds along the path $\gamma(\cdot)$,   then it also holds at each point of a connected  open neighborhood $\mathcal{O}$ of the path $\gamma(\cdot)$ in $M$. Applying Rashevsky-Chow theorem to the restriction of the control system (\ref{cls}) to
 $\mathcal{O}$ we get the needed steering result.

In the case of continual ensembles it turns out though   - and this is the second difficulty -  that for the vector fields $F$, which are  lifts to  $\mathcal{E}_\Theta(M)$ of the vector fields $f \in \mbox{Vect }M$ ,  the (approximate)  bracket generating property \eqref{bragen}
 can not hold at each $\gamma \in \mathcal{E}_\Theta(M)$  and may cease  to hold  even  $C^0$-locally.
Thus the argument just provided fails:    condition  (\ref{bragen}) may  hold along the path $p(\cdot)$ and
cease to hold in a neighborhood of the path.

For example the space $\mathcal{E}=\mathcal{E}_\Theta(\mathbb{R}^n)$ of  ensembles of points in $\mathbb{R}^n$,  parameterized by a compact
$\Theta$,
 is isomorphic to the
Banach space $C^0(\Theta , \mathbb{R}^n)$.  Its tangent spaces are all isomorphic  to $C^0(\Theta ,\mathbb{R}^n)$. If $\Theta$ is not finite ($\sharp \Theta = \infty$) then in any $C^0$-neighborhood
     of an ensemble $\hat \gamma(\cdot) \in C^0(\Theta ,\mathbb{R}^n)$ one can find an ensemble  $\gamma(\cdot) \in C^0(\Theta ,\mathbb{R}^n)$, which is constant on an open subset of $\Theta$. Then  $\{Y(\gamma(\theta))|  Y \in \mbox{Vect} M\}$  is not  dense in $T_\gamma \mathcal{E}=T_\gamma C^0(\Theta ,\mathbb{R}^n)$ and hence condition \eqref{bragen} can not hold at $\gamma(\cdot)$. There may certainly occur other types of singularities.

The same remains true if the topology, in which the target is approximated (and hence the topology of $\mathcal{E}$) is weakened.

We end up with two remarks concerning the formulation of Theorem \ref{trc}.

The criterion for approximate steering, provided by the Theorem has meaningful analogue also in the case of singletons.

\begin{prop}[bracket approximating property and approximate steering for singletons]\label{baac} Let   $x_\alpha , x_\omega \in M$ and $\gamma(t), \ t \in [0,T]$
 be a   continuously differentiable  path, which joins  $x_\alpha$ and $x_\omega$. If the Lie bracket approximating property holds at each point $\gamma(t), \ t \in [0,T]$,  then $x_\alpha$ can be approximately steered to $x_\omega$ by an admissible trajectory of $(\ref{cls})$. $\square$
\end{prop}

 Recall that the Lie bracket approximating condition includes the assumption   of
 Lipschitz equicontinuity
of the approximating vector fields from $\mbox{Lie}\{f\}$.  The following example  illustrates     importance of  this  assumption.

Consider a control system (\ref{cls}) in $\mathbb{R}^2=\{(x_1,x_2)\}$,
   such that the orbits of \eqref{cls} or, the same, of the Lie algebra $\mbox{Lie}\{f\}$   are the lower and the upper open half-planes of $\mathbb{R}^2$
   together with  the straight-line $x_2=0$.
   The points  $x_\alpha=(-1,-1)$ and $x_\omega=(1,1)$ belonging to different orbits,  can not be steered approximately one to another.
On the other side if we join these  points  by the curve $\gamma(t)=(t,t^3), \ t \in [-1,1]$, then it is  immediate to check, that
$\dot{\gamma}(t) \in \mbox{Lie}\{f\}(\gamma(t))$ for each $t$, but the condition of Lipshitz equicontinuity is
not  fulfilled.
There are curves $\gamma^\delta(\cdot)$ arbitrarily close to $\gamma(\cdot)$ in $C^0$ metric,
which intersect the line $x_2=0$ transversally and hence do not satisfy
the condition
$\dot{\gamma}^\delta (t) \in \mbox{Lie}\{f\}(\gamma^\delta(t))$.

\section{Proof of Theorem \ref{generic}.}
\label{proof1}
We provide a proof for couples of vector fields ($s=2$); general case is treated similarly.
 It suffices to establish for fixed $N$  existence of a residual subset
 $\mathcal{G} \subset \mbox{\rm Vect} M \times \mbox{\rm Vect} M$ such that
 for each couple $(X,Y) \in \mathcal{G}$ the couple of $N$-folds of the vector fields $(X^N, Y^N)$ is bracket generating on $M^{(N)}$.
 Let  $\dim M=n$.

The proof is based on  application of  J.Mather's multi-jet transversality theorem (\cite{GG}).
 Consider
the couples of vector fields $(X,Y)$ on $M$ as  $C^k$-smooth sections of  the fibre bundle $\pi:TM\times_M TM \to M$. Consider the set
$J_k(TM\times_M TM)$ of $k$-jets of the couples of vector fields and the projection  $\pi_k$  of $J_k(TM\times_M TM)$ to $M$. One can define in obvious way
for $N \geq 1$ the projection
$\pi_k^N:J_k(TM \times_{M} TM)^N \to M^N$ and introduce  the set $J_k^{(N)}(TM \times_{M} TM)^N=(\pi_k^N)^{-1}(M^{(N)})$,
 which is $N$-fold $k$-jet (or  multi-jet) bundle
for the  couples of vector fields.

In other words $N$-fold of a vector field $X \in \mbox{Vect}M$ is a vector field $\underbrace{(X, \ldots , X)}_N \in \mbox{Vect}M^{(N)}$.  For a couple $(X,Y) \in \mbox{Vect} M \times \mbox{Vect} M$ of  vector fields the multi-jet  $J_k^{(N)}(X,Y): M^{(N)} \to J_k^{(N)}(\mbox{Vect} M \times \mbox{Vect} M)$  can be represented as
\begin{eqnarray*}
  \forall (x_1, \ldots , x_N) \in M^{(N)}:  \\
  J_k^{(N)}(X,Y)(x_1, \ldots x_N)=\left(J_k(X,Y)(x_1), \ldots , J_k(X,Y)(x_N) \right).
\end{eqnarray*}


\begin{prop}[multi-jet transversality theorem for the couples of vector fields; \cite{GG}]\label{multijet}
 Let $S$ be a submanifold  of the space of $k$-multijets ($N$ fold $k$-jets) $J_k^{(N)}(TM \times_{M} TM)^N$.  Then for sufficiently large $\ell$ the set of couples of the vector fields
 \[T_S=\{(X,Y) \in \mbox{\rm Vect} M \times \mbox{\rm Vect} M| \ J_k^N(X,Y) \transv S\} \]
 is a residual subset of  $\mbox{\rm Vect} M \times \mbox{\rm Vect} M$ in  Whitney $C^\ell$-topology ($\transv$ stays for transversality of a map to a manifold). $\Box$
\end{prop}

Coming back to the proof of Theorem \ref{generic}, note that
the set $\mathcal{R}$ of the couples $(X,Y)$ of vector fields, such that at each $x \in M$ either $X(x) \neq 0$, or $Y(x) \neq 0$,
 is open and dense in   $\mbox{Vect}M \times \mbox{Vect}M$.
We will seek $\mathcal{G}$ as a subset of  $\mathcal{R}$.

For each couple $(X,Y) \in \mathcal{R}$,
and each point $\bar x=(x_1, \ldots , x_N) \in M^{(N)}$  we introduce  the two $nN \times 2nN$-matrices:
\begin{eqnarray*}
  V(\bar x)=\left(
      \begin{array}{cccc}
        Y(x_1) & \ad X Y (x_1) & \cdots  & \ad^{2nN-1} X Y(x_1) \\
        \vdots & \vdots  & \vdots  & \vdots  \\
        Y(x_N) & \ad X Y (x_N) & \cdots  & \ad^{2nN-1} X Y(x_N) \\
      \end{array}
    \right), \\
   W(\bar x)=\left(
      \begin{array}{cccc}
        X(x_1) & \ad^2 Y X (x_1) & \cdots  & \ad^{2nN} Y X(x_1) \\
        \vdots & \vdots  & \vdots  & \vdots  \\
        X(x_N) & \ad^2 Y X (x_N) & \cdots  & \ad^{2nN} Y X(x_N) \\
      \end{array}
    \right).
 \end{eqnarray*}
(Note that $W(\bar x)$ lacks the column constituted by $\ad Y X (x_j)$ which coincides, up to a sign,  with the second column in $V(\bar x)$).

For  $(X,Y) \in \mathcal{R}, \ \bar x=(x_1, \ldots , x_N) \in M^{(N)} $ and  each $ x_i, \ i=1, \ldots N,$ at least one of the vectors $X( x_i),Y( x_i)$ is non null. We can choose local coordinates
$\xi_{ij}, \ i=1, \ldots N; \ j=1, \dots n$  in a
neighborhood $U=U_1 \times \cdots \times  U_N$ of $\bar x=(x_1, \ldots , x_N) \in M^{(N)}$ in such a way that  in each $U_i, \ i=1, \ldots ,N$ either
$X$ or $Y$ becomes the non null constant  vector field:   $X=\partial/\partial \xi_{i1}$ or $Y=\partial/\partial \xi_{i1}$.
Then  for each $i=1, \ldots ,N,$ either $\ad^{k} X Y|_{x_i}$ or  $\ad^{k} Y X|_{x_i}$ equal respectively  to $\left.\frac{\partial^k Y}{\partial \xi_{i1}^k}\right|_{x_i}$ or
$\left.\frac{\partial^kX}{\partial \xi_{i1}^k}\right|_{x_i}$.

We call {\it significant} those elements of the $(Nn \times 2Nn)$-matrices $V(\bar x)$,  $W(\bar x)$
and of the corresponding $(Nn \times 4Nn)$-matrix $(V(\bar x)|W(\bar x))$, which are the  components of
$\frac{\partial^kY}{\partial \xi_{i1}^k}$  and of
$\frac{\partial^k X}{\partial \xi_{i1}^k}$.
For each $j=1, \ldots ,Nn$  either $j$-th row of $V(\bar x)$ or $j$-th row of $W(\bar x)$  consists of significant elements. The elements of these matrices are polynomials in the components of the multi-jets $J^{2nN}X(\bar x),J^{2nN}Y(\bar x)$. Significant elements are polynomials of degree $1$, distinct significant elements correspond to different polynomials,  nonsignificant elements correspond to polynomials of degrees $>1$.
Elements of  different rows of the matrices  differ.

If $(X,Y) \in \mathcal{R}$  and $(X^N,Y^N)$ lacks the bracket generating property at some $\bar x=(x_1, \ldots , x_N)$, then  the rank $r$
of the $(Nn \times 4Nn)$-matrix $(V|W)(\bar x)$ is incomplete:
  $r< nN$.

The (stratified) manifold of $(Nn \times 4Nn)$-matrices of rank $r <nN$ is (locally) defined by rational relations, which express elements of some
$(Nn-r)\times (4Nn-r)$ minor via other elements of the matrix.

As long  as $4Nn-r \geq 3Nn+1$,  then each row of the minor contains  $s \geq 3Nn+1 - 2Nn >Nn$ significant
elements. The corresponding  relations express  $s$
distinct components  of $2N$-th multi-jet of $(X,Y)$  via other components of the multi-jet. Hence $2N$-multi-jets of the couples $(X,Y)$, for which $(X^N,Y^N)$
lack bracket generating property,
must belong to an algebraic  manifold $S$ of codimension $s> Nn$ in $J_k^N(TM \times_M TM$.

Consider the set $T_S$ of the couples $(X,Y) \in \mathcal{R} \subset \mbox{\rm Vect} M \times \mbox{\rm Vect} M$,
for which  $J_{2nN}^N(X,Y):M^{(N)} \to J_{2nN}^N(\mbox{Vect} M \times \mbox{Vect} M)$ is transversal to $S$. According to the multijet transversality theorem
(Proposition \ref{multijet})  $T_S$  is residual in $\mbox{\rm Vect} M \times \mbox{\rm Vect} M$ in
Whitney $C^\ell$-topology for sufficiently large $\ell$.
As far as $\dim M^{(N)} =Nn<s=\mbox{\rm codim } S$, the transversality
can take place only if,  for each $\bar x \in M^{(N)}$,   $J_{2nN}^N(X,Y)|_{\bar x} \not\in S$.
Hence for each couple $(X,Y)$ from the residual  subset  $T_S$,
 the  couples of $N$-folds $(X^N,Y^N)$ are bracket generating at each  point of $M^{(N)}$.

\section{Proof of Theorem  \ref{liex_flows}}\label{proof3}
\subsection{Variational formula}

We start with     nonlinear version of 'variation of constants' formula, which will be employed in the next subsection.

 Let $f_t(x)$ be a time-variant and $g(x)$ a  time-invariant vector fields on $M$. We assume both vector fields to be $C^\infty$-smooth and   Lipschitz on $M$.
Let $\chro_0^t f_\tau d\tau$ denote the flow generated by the time-variant vector field $f_t$ (see \cite{AGam77,AgSach} for the notation),
and  $e^{tg}$  stays for  the flow, generated by the time-invariant vector field $g$.

 \begin{lem}[\cite{AgSach}]\label{varf}
   Let $f_\tau(x),g(x)$ be $C^\infty$-smooth  in $x$, $f_\tau$ integrable in $\tau$. Let $U(t)$ be a Lipschitzian function on $[0,T], \ U(0)=0$.
   The flow $P_t=\chro_0^t \left( f_\tau(x)+g(x)\dot{U}(\tau) \right) d\tau$, generated by
    the differential equation
\begin{equation}\label{caf}
  \dot{x}=f_t(x)+g(x)\dot{U}(t),
\end{equation}
   can be represented as a composition of flows
  \begin{equation}\label{vf}
    \chro_0^t \left(f_\tau(x)+g(x)\dot{U}(\tau)\right)d\tau=\chro_0^t \left(e^{-U(\tau) g}\right)_*f_\tau d\tau \circ e^{U(t)g}.    \qed \end{equation}
   \end{lem}
      At the right-hand side of \eqref{vf} $\left(e^{-U(\tau) g}\right)_*$ is the differential of the diffeomorphism  $e^{-U(t) g}=\left(e^{U(t) g}\right)^{-1},$ where $e^{U(t) g}$ is  the evaluation at time-instant $U(t)$  of the flow,  generated by the time-invariant vector field $g(x)$.

We omit at this point  the questions of  completeness of the vector fields involved into \eqref{caf},\eqref{vf}, assuming that the formula \eqref{vf} is valid, whenever the flows, involved in it, exist on the specified intervals.

 For each vector field $Z \in \mbox{Vect }M$ the  operator $\ad_{Z}$,
acts on the space of vector fields:
$\ad_Z Z_1=[Z,Z_1]$ -  the Lie bracket of $Z$ and $Z_1$.   The operator exponential $e^{U \ad Z}$ is defined formally:
 $e^{U\ad_Z}=\sum_{j=0}^\infty \frac{U^j(\ad_Z)^j}{j!}$.
For $C^\infty$-smooth  vector fields $Z,Z_1$ the expansion is known (see \cite{AGam77},\cite{AgSach}) to provide  asymptotic  representation for $\left(e^{-U(\tau) g}\right)_*$: for each $s \geq 0$ and a
 compact $K \subset M$ there exists a compact neighborhood $K'$ of $K$ and $c>0$ such that
 \begin{eqnarray*}
   \left\|\left(\left(e^{-U(\tau) g}\right)_*-I-\sum_{j=0}^{N-1} \frac{(U(\tau))^j}{j!}\ad^j_g\right)Z_1\right\|_{s,K} \leq  \nonumber \\
  \leq   ce^{c|U(\tau)|\|g\|_{s+1,K'}}\frac{\left(|U(\tau)|\|g\|_{s+N,K'}\right)^N}{N!}\|Z_1\|_{s+N,K'}
 \end{eqnarray*}
(see \cite{AGam77} for the details).
We  employ   the asymptotic  formulae  for $N=1,2$ and small $|U|$:
\begin{eqnarray}\label{2expad}
  \left\|\left(\left(e^{-U(\tau) g}\right)_*-I\right)Z_1\right\|_{s,K}=O(|U(\tau)|)\|Z_1\|_{s+1,K'}, \\
  \left\|\left(\left(e^{-U(\tau) g}\right)_*-I-U(\tau)\ad_g \right)Z_1\right\|_{s,K}=o(|U(\tau)|)\|Z_1\|_{s+2,K'},
\label{22expad}
\end{eqnarray}
 as $|U| \to 0$.

 We introduce  at this point fast-oscillating controls by  choosing  $1$-periodic Lipschitz  function $V(t)$ with $V(0)=0$,  the scaling parameters $\beta > \alpha >0$ and defining for  $\eps>0$:
 $V(t;\alpha , \beta , \eps )=\eps^\alpha V\left(t/\eps^\beta\right)$.  We introduce
  controls    $$u_\eps(t)=\frac{dV(t;\alpha , \beta , \eps )}{dt}=\eps^{\alpha-\beta} \dot{V}\left(t/\eps^\beta\right),$$
  which are high-gain and   fast-oscillating for small $\eps>0$.

  Fore a more general control
  \begin{equation}\label{ueps}
   u_\eps(t)=w(t)\eps^{\alpha - \beta}\dot{V}\left(t/\eps^\beta\right),
  \end{equation}
   where $w(\cdot)$
  is a  Lipschitz function, the primitive of    $ u_\eps(t)$  equals
  \begin{equation}\label{uprim}
    U_\eps (t) =\eps^{\alpha}\left(w(t)V\left(t/\eps^\beta\right)-\int_0^t V\left(\tau/\eps^\beta\right)\dot{w}(\tau) d\tau \right)=\eps^{\alpha}\hat U_\eps(t),
  \end{equation}
  and $\hat U_\eps (t)=O(1)$   as $\eps \to +0$ uniformly for $t$ in a compact interval.

Substituting    $U(t)=U_\eps(t)$, defined by \eqref{uprim}, into  \eqref{vf} we get
\begin{eqnarray}\label{exuw}
  \chro_0^t \left(f_\tau(x)+g(x)\eps^{\alpha-\beta} w(\tau)\dot{V}\left(\frac{\tau}{\eps^\beta}\right)\right)d\tau=  \\   =\chro_0^t \left(e^{-\eps^{\alpha} \hat U_\eps(\tau)g}\right)_*f_\tau d\tau \circ e^{\eps^{\alpha}  \hat U_\eps (t)g}. \nonumber
\end{eqnarray}


Expanding the exponentials  at the right-hand side of the   equality   according to  formula \eqref{2expad}    we get for
the   control $ u_\eps(t)$, defined by \eqref{ueps}:
\begin{align}\label{fosc}
  \chro_0^t \left(f_\tau(x)+g(x)u_\eps(\tau)\right)d\tau=  \nonumber \\
  =\chro_0^t \left( f_\tau(x) +O(\eps^{\alpha})\right)d\tau \circ  \left(I+O(\eps^{\alpha})\right).
\end{align}

By classic theorems on continuous dependence of trajectories on the right-hand  side we conclude
that the  flow $\chro_0^t \left(f_\tau(x)+g(x)u_\eps(\tau)\right)d\tau$ with $u_\eps(t)$, defined by \eqref{ueps}, tends to $\chro_0^t f_\tau(x)d\tau$, as $\eps \to 0$,
uniformly in $t$ on compact intervals. Therefore the effect of the fast-oscillating control \eqref{ueps}  tends to zero  as $\eps \to 0$
with respect to  any   of the seminorms $\|\cdot\|_{r,K}$:
\[\left\|\chro_0^t \left(f_\tau(x)+g(x)u_\eps(\tau)\right)d\tau
  -\chro_0^t f_\tau(x)d\tau  \right\|_{r,K}  \rightrightarrows  0\] for all $r \geq 0$,   compact $K$ and uniformly for $t \in [0,T].  \hfill \qed$

\subsection{Lie extension for flows}

Coming back to the proof of  Theorem \ref{liex_flows} we first note that the conclusion  can be arrived at  by induction, with the step of induction, represented by
the following

\begin{lem}
The conclusion of the theorem holds for the controlled  system
  \[\frac{d}{d t} x(t) =\sum_{j=1}^kX^j(x)u_j(t)+    X(x) u(t)+Y(x) v(t),\]
and its Lie extension
  \[\frac{d}{d t} x(t) =\sum_{j=1}^kX^j(x)u^e_j(t)+ X(x) u^e(t)+Y(x) v^e(t)+[X , Y](x)w^e(t).\]
\end{lem}

The proof, provided below,  shows  that one can
leave out,  without loss of generality,  the summed addends
$\sum_{j=1}^kX^k(x)u_k(t)$,  $\sum_{j=1}^kX^k(x)u^e_k(t)$ at the right-hand side of the systems.  It suffices to prove the  result
 for the  $2$-input system
\begin{equation}\label{eq_6}
  \frac{d}{d t} x(t) =    X(x) u(t)+Y(x) v(t),
\end{equation}
and its $3$-input Lie extension
\begin{equation}\label{eq_7}
  \frac{d}{d t} x(t) =X(x) u^e(t)+Y(x) v^e(t)+[X , Y](x)w^e(t).
\end{equation}
One can assume,  without loss of  generality,  $w^e(t)$  to be  smooth,  as far as smooth functions
  are dense  in $L_1$-metric in the space of bounded measurable functions.
  Hence by classical results on continuous dependence with respect to right-hand sides,
   the  flows, generated by measurable  controls,
  can be approximated by    flows, generated by smooth  controls.

To construct  the controls $u(t), v(t)$  from $u^e(t), v^e(t), w^e(t)$   we take
\begin{equation}\label{econ}
  u(t)=u_\eps(t)=u^e(t)+\eps \dot{U}_\eps(t), \ v(t)=v_\eps(t)=v^e(t)+\eps^{-1}\hat{v}_\eps(t),
\end{equation}
where $\eps$ is the  parameter of  approximation and the  functions $U_\eps(t)$ and $\hat v_\eps (t)$ will be specified in a moment.

Feeding  controls  \eqref{econ} into   system \eqref{eq_6} we get
\begin{equation}\label{eq_8}
  \frac{d}{d t} x(t)=\underbrace{X(x)u^e(t)+Y(x)\left(v^e(t)+\eps^{-1}\hat{v}_\eps(t)\right)}_{f_t}+\underbrace{X(x)}_g \eps  \dot{U}_\eps(t).
\end{equation}
Applying formula \eqref{vf} to the flow, generated by  \eqref{eq_8}, we represent it as a
composition
\begin{align}\label{crex}
\chro_0^t X(x)u^e(t)+\left(e^{-\eps U_\eps(t)  X}\right)_*Y(x)\left(v^e(t)+\eps^{-1}\hat{v}_\eps(t)\right)dt \ \circ  \nonumber \\  \circ  \ e^{\eps U_\eps(t) X(x)}.
\end{align}

We wish the latter flow   to   approximate (for sufficiently small $\eps >0$) the flow, generated by   \eqref{eq_7}.
To achieve this we  choose the functions
\begin{equation}\label{uveps}
  U_\eps(t)=2\sin(t/\eps^2)w^e(t), \ \hat{v}_\eps(t)=\sin(t/\eps^2).
\end{equation}

Approximating the operator exponential $e^{\eps U_\eps(t) \ad_X} $     by formula \eqref{22expad}  we
  transform \eqref{crex} into
\begin{align}\label{xeps}
 \!\! \chro_0^t(X (x)u^e(t)+&  Y (x)v^e(t)+  [X,Y](x)U_\eps(t)\hat{v}_\eps(t)+\\
 +Y(x) \eps^{-1}\hat{v}_\eps(t)+&O(\eps))dt \circ (I+O(\eps)) ,  \nonumber
\end{align}
where all $O(\eps)$ are uniform in $t \in [0,T]$.

From     \eqref{uveps}
\[ U_\eps(t)\hat{v}_\eps(t)=w^e(t)-w^e(t)\cos(2t/\eps^2), \]
and  \eqref{xeps} takes form
\begin{align}\label{xxeps}
  \chro_0^t\!\! \left(X (x)u^e(t)+ Y (x)v^e(t) +[X,Y](x)w^e(t)+   Y(x) \eps^{-1}\sin(t/\eps^2) - \right. \\  -    \left. [X,Y](x)w^e(t)\cos(2t/\eps^2)+O(\eps)\right) dt  \circ (I+O(\eps)) . \nonumber
\end{align}


Processing fast oscillating terms $ Y(x)\eps^{-1}\sin(t/\eps^2), \   [X,Y]w^e(t)\cos(2t/\eps^2)$  according to formula \eqref{exuw}
we bring  the flow \eqref{xxeps} to the form
\begin{align*}
  \chro_0^t\!\! \left(X (x)u^e(\tau)+ Y (x)v^e(\tau) +[X,Y](x)w^e(\tau)+O(\eps)\right) d\tau  \ \circ \\
  \circ  \ (I+O(\eps)) ,
\end{align*}
 wherefrom  one concludes  for $u_\eps(t),v_\eps(t)$,  defined by formulae   \eqref{econ}-\eqref{uveps}, the  convergence of the flows: for each $r \geq 0$ and compact $K$
   \begin{eqnarray*}
      \left\|\chro_0^t\!\! \left(X (x)u^e(\tau)+ Y (x)v^e(\tau) +[X,Y](x)w^e(\tau) \right) d\tau - \right. \\
     -  \left. \chro_0^t \left(X(x) u_\eps(\tau) +Y(x) v_\eps(\tau)\right)d\tau\right\|_{r,K}=O(\eps)
    \end{eqnarray*}
 as $
     \eps \to 0$.

\section{Proof of Theorem \ref{trc}.}
\label{proof2}

\subsection{Steering ensembles of  points by an extended control}

\begin{prop}\label{gloec}
Under the assumptions of Theorem  \ref{trc},   for each $ \eps >0$
there exists a finite set $B$ (depending on $\eps$) of the multiindices  $\beta=(\beta_1, \ldots , \beta_N)$ and an  extended differential equation (\ref{ext_ens})
together with an extended control $(v_\beta(t))_{\beta \in B}, \ t \in [0,T]$ such that the flow, generated  by (\ref{ext_ens}) and
the control steers,  in time $T$,  the initial ensemble $\alpha(\theta)$ to the ensemble $x(T;\theta)$, for which
$\sup_{\theta \in \Theta} d\left(x(T;\theta), \omega(\theta)\right) < \eps. \  \qed$
\end{prop}



Consider the diffeotopy $\gamma_t(\theta)=P_t (\alpha(\theta))$, along which
 Lie bracket $C^0$-approximating condition   holds. Let $Y_t(x)$  be the time-variant vector field,  which generates the diffeotopy and $\Gamma$ its image.
We start with  the following  technical Lemma.

\begin{lem}\label{contappr}
Let  assumptions of  Theorem~\ref{trc}  hold. Then there exists $\lambda >0$ and  compact neighborhood $W_\Gamma \supset \Gamma$,  such that
for each $\eps>0$  there exists a
 finite set of multi-indices $B$ together with   continuous
functions $\left(v_\beta(t)\right), \ \beta \in B$ such that  $X_t(x)=\sum_{\beta \in B}v_\beta(t)X^\beta(x)$  satisfies:
\begin{equation}\label{YvX}
  \|X_t(x)\|_{1,W_\Gamma} < \lambda ,  \ \left\|Y_t(\gamma_t(\theta)) -X_t(\gamma_t(\theta)) \right\|_{C^0(\Theta)}< \eps . \ \Box
\end{equation}
\end{lem}

{\it Proof of Lemma~\ref{contappr}.}  According to the Lie bracket $C^0$-approximating assumption along the diffeotopy there exists
 $\lambda>0$ and for each $t \in [0,T]$ and each $\eps>0$ a
finite set $B_t$ of multi-indices and the coefficients $c_\beta(t), \ \beta \in B_t,$ such that
\begin{eqnarray}\label{appr_Y}
    \left\|\sum_{\beta \in B_t}c_\beta(t)X^\beta(x)\right\|_{1,W_\Gamma} < \lambda ,  \nonumber \\ \left\|Y_t(\gamma_t(\theta)) -\sum_{\beta \in B_t}c_\beta(t)X^\beta(\gamma_t(\theta)) \right\|_{C^0(\Theta)}< \eps  .
\end{eqnarray}

As far as $Y_t(\gamma_t(\theta))$ and $X^\beta(\gamma_t(\theta))$ vary continuously with $t$,
the estimate
\begin{equation*}
\left\|Y_\tau (\gamma_\tau (\theta)) -\sum_{\beta \in B_t}c_\beta(t)X^\beta(\gamma_\tau (\theta)) \right\|_{C^0(\Theta)}< \eps   \end{equation*}
is valid for $\tau \in  \mathcal{O}_t$ - a neighborhood of $t$. The family $\mathcal{O}_t \ (t \in [0,T])$ defines an open covering of $[0,T]$,
from which we choose finite subcovering $\mathcal{O}_i=\mathcal{O}_{t_i}, \ i=1, \ldots , N$.
Putting  $B=\bigcup_{i=1}^N B_{t_i}$ we define  $c_{i\beta}=c_\beta(t_i), \ \forall i=1, \ldots , N, \   \forall \beta \in B_i$.

Choose  a smooth partition of unity $\{\mu_i(t)\}$ subject to the covering $\{\mathcal{O}_i\}$.
Put for each $\beta \in B, \ v_\beta(t)=\sum_{i=1}^N\mu_i(t)c_{i\beta}$;
it is immediate  to see that $v_\beta(t)$ are continuous.   For
\begin{equation}\label{vec_xt}
  X_t(x)=\sum_{\beta \in B}v_\beta (t)X^\beta(x)
\end{equation}
we conclude
\begin{eqnarray*}
 \forall \theta \in \Theta : \  \left\|Y_t(\gamma_t(\theta)) - X_t(\gamma_t(\theta))\right\|= \\
  =\left\| \sum_{i=1}^N\mu_i(t)Y_t(\gamma_t(\theta)) -
  \sum_{i=1}^N \sum_{\beta \in B}\mu_i(t)c_{i\beta}X^\beta (\gamma_t(\theta))\right\| &\leq& \\
  \leq \sum_{i=1}^N \mu_i(t) \left\| Y_t(\gamma_t(\theta)) - \sum_{\beta \in B_i}c_{i\beta}X^\beta (\gamma_t(\theta))\right\| \leq \eps \sum_{i=1}^N \mu_i(t)&=&\eps .
\end{eqnarray*}
The first of the estimates \eqref{YvX} is proved similarly. $\qed$




Coming back to the proof of  Proposition \ref{gloec} we  consider the evolution of the ensemble $\alpha(\theta)$
under the action of the flow generated by  the vector field $X_t$, defined by \eqref{vec_xt}.


We   estimate
\begin{eqnarray*}\|x(t; \theta)-\gamma_t(\theta))\|=\left\|\int_0^t\left(X_\tau(x(\tau;\theta), v(\tau)))-Y_\tau(\gamma_\tau(\theta))\right)d\tau  \right\| \leq \\
\leq  \int_0^t\left\|X_\tau(x(\tau;\theta))-X_\tau(\gamma_\tau(\theta))\right\|d\tau  +  \int_0^t\left\|X_\tau(\gamma_\tau(\theta))-Y_\tau(\gamma_\tau(\theta))\right\|d\tau  .
\end{eqnarray*}

By virtue of (\ref{appr_Y}) we obtain (whenever $x(t;\theta) \in W_\Gamma$):
\[\|x(t;\theta)-\gamma_t(\theta)\|   \leq \lambda\int_0^t  \left\|x(\tau;\theta)-\gamma_\tau(\theta)\right\| d\tau + \eps t ,\]
and by Gronwall lemma
\begin{equation}\label{xgamma}
\|x(t;\theta)-\gamma_t(\theta)\|  \leq \eps \frac{\left(e^{\lambda t}-1\right)}{\lambda}.
\end{equation}

We should take $\eps$ sufficiently small, so that \eqref{xgamma} guarantees that $x(t;\theta)$ does not leave the neighborhood $W_\Gamma$, defined by Lemma \ref{contappr}.  Then
\[\|x(T;\theta)-\omega(\theta)\|  \leq \eps \frac{\left(e^{\lambda T}-1\right)}{\lambda}  \]
  and the claim of Proposition \ref{gloec} follows. $\qed$

Theorem \ref{trc} follows readily from Propositions  \ref{gloec} and    Corollary \ref{liex_ens}.

\bibliographystyle{model1-num-names}



\end{document}